\documentclass[12pt]{article}
 
\usepackage{txfonts}
\usepackage[american]{babel}

\def\II{{\mathbb I}} 
\def\RR{{\mathbb R}} 
\def\CC{{\mathbb C}}

\def\tr{\mathrm{ tr\,}} 
\def\Tr{\mathrm{ Tr\,}}

\def\vol{\mathrm{ vol\,}}

\def\diag{\mathrm{diag\,}} 
\def\Spin{\mathrm{Spin}}

\def\be{\begin{equation}} 
\def\ee{\end{equation}} 
\def\bea{\begin{eqnarray}} 
\def\eea{\end{eqnarray}} 
\def\bed{\begin{definition}{\ }}
\def\eed{\end{definition}}
\def\bd{\begin{description}}
\def\ed{\end{description}}
\def\bc{\begin{center}}
\def\ec{\end{center}}

\newtheorem{theorem}{Theorem}

\newtheorem{proposition}{Proposition}

\newtheorem{definition}{Definition}

\begin{document}

\begin{titlepage}

\null
\vskip-2cm 
\hfill
\begin{minipage}{7cm}
\par\hrulefill\par\vskip-4truemm\par\hrulefill
\par\vskip2mm\par
{{\large\sc New Mexico Tech 
{\rm (July 2007)}}}
\par\hrulefill\par\vskip-4truemm\par\hrulefill
\end{minipage}
\bigskip 
\bigskip

\hfill
\begin{minipage}{7cm}
to be published in the special volume of {\it SIGMA},
Proceedings of the {\it 2007 Midwest Geometry Conference}
in Honor of Thomas P. Branson, Iowa City, IA, May 18-20, 
2007
\end{minipage}

\par
\bigskip
\bigskip

\centerline{\LARGE\bf Heat Kernel Asymptotics on Homogeneous Bundles}
\bigskip
\bigskip
\bigskip
\centerline{\Large\bf\textbf{Ivan G. Avramidi}}
\bigskip
\centerline{\it Department of Mathematics}
\centerline{\it New Mexico Institute of Mining and Technology}
\centerline{\it Socorro, NM 87801, USA}
\bigskip

\centerline{\sc 
Dedicated to the Memory of Thomas P. Branson}
\medskip 

\begin{abstract}
We consider Laplacians acting on sections of homogeneous vector bundles
over symmetric spaces.  By using an integral representation of the heat
semi-group  we find a formal solution for the heat kernel diagonal that
gives a generating function for the whole sequence of heat invariants.
We argue that the obtained formal solution correctly
reproduces the exact heat kernel diagonal after a suitable
regularization and analytical continuation.
\end{abstract} 

\vfill



\end{titlepage}


\section{Introduction}
\setcounter{equation}0

The heat kernel is one of the most powerful tools in mathematical
physics and  geometric analysis. Of special importance is the 
short-time  asymptotic expansion
of the trace of the heat kernel.
The coefficients of this asymptotic
expansion, called the  heat invariants, are extensively used in
geometric analysis, in particular, in spectral geometry and index
theorems proofs. There has been a tremendous progress in the explicit
calculation  of spectral asymptotics in the last thirty years
\cite{gilkey75,avramidi89,avramidi90,avramidi91,vandeven98,
yajima04}
(see also the reviews \cite{gilkey95,avramidi99,
avramidi02,vassilevich03,kirsten01}).  
A further progress in the study of spectral
asymptotics can be only achieved by restricting oneself to operators and
manifolds with high level of symmetry, in particular, homogeneous
spaces, which enables one to employ powerful algebraic methods.  
It is well known that heat invariants are  determined essentially by
local geometry. They are polynomial invariants in the curvature with
universal constants that do not depend on the global properties of the
manifold. It is this universal structure that we are
interested in this paper.   Our goal is to compute the heat kernel
asymptotics of the Laplacian acting on homogeneous vector bundles over
symmetric spaces. 

\section{Geometry of Symmetric Spaces}
\setcounter{equation}{0}

\subsection{Twisted Spin-Tensor Bundles}

In this section we introduce basic concepts and fix notation. Let 
$(M,g)$ be an $n$-dimensional Riemannian manifold without boundary. We
assume that it is complete simply connected orientable and spin. We
denote the local coordinates on $M$ by $x^\mu$, with Greek indices
running over $1,\dots, n$. Let $e_{a}{}^\mu$ be a local orthonormal
frame defining a basis for the tangent space $T_xM$. We denote the frame
indices by low case Latin indices from the beginning of the alphabet,
which also run over $1,\dots,n$. The frame indices are raised and
lowered by the metric $\delta_{ab}$.  Let $e^{a}{}_\mu$ be the matrix
inverse to $e_{a}{}^\mu$, defining the dual basis in the cotangent space
$T_x^*M$.  As usual, the orthonormal frame, $e^a{}_\mu$ and $e_a{}^\mu$, 
will be used to transform the coordinate (Greek) indices to the
orthonormal (Latin) indices.  The covariant derivative  along the frame
vectors is  defined by $\nabla_a=e_a{}^\mu\nabla_\mu$. For example, with
our notation, $\nabla_a\nabla_b T_{cd} =e_a{}^\mu e_b{}^\nu e_c{}^\alpha
e_d{}^\beta\nabla_\mu\nabla_\nu  T_{\alpha\beta}$.

Let $\mathcal{T}$ be a spin-tensor bundle realizing a representation
$\Sigma$  of the spin group \mbox{$\Spin(n)$} with the fiber $\Lambda$
and let $\Sigma_{ab}$ be the generators of the its Lie algebra. The spin
connection naturally induces a connection  on the bundle $\mathcal{T}$.
The curvature of this connection is $\frac{1}{2}R^{ab}{}\Sigma_{ab}$,
where $R^{ab}{}$ is the curvature 2-form of the spin connection.

In the present paper we will further  assume that $M$ is a {\it locally
symmetric space} with a Riemannian metric with the parallel curvature,
that is,
$\nabla_\mu R_{\alpha\beta\gamma\delta}=0$,
which means, in particular, that the Riemann curvature tensor
satisfies the 
integrability constraints
\be
R^{fg}{}_{ e a}R^e{}_{b cd} 
-R^{fg}{}_{ e b}R^e{}_{a cd} 
+ R^{fg}{}_{ ec}R^e{}_{d ab}
-R^{fg}{}_{ ed}R^e{}_{c ab}
= 0\,.
\label{312}
\ee


Let $G_{YM}$ be a  compact Lie group (called a gauge group). It
naturally defines the principal fiber bundle over the manifold $M$  with
the structure group $G_{YM}$. We consider a representation of the
structure group $G_{YM}$ and the associated vector bundle through this
representation with the same structure group $G_{YM}$ whose typical
fiber is a $k$-dimensional  vector space $W$. Then for any spin-tensor
bundle $\mathcal{T}$ we define the twisted spin-tensor bundle
$\mathcal{V}$ via the twisted product of the bundles $\mathcal{W}$ and
$\mathcal{T}$. The fiber of the bundle $\mathcal{V}$ is
$V=\Lambda\otimes W$ so that the sections of the bundle $\mathcal{V}$
are represented locally by $k$-tuples of spin-tensors.

A connection on the bundle $\mathcal{W}$ (called Yang-Mills or gauge
connection) taking values in the Lie algebra $\mathcal{G}_{YM}$ of the
gauge group $G_{YM}$ naturally defines the total connection on the
bundle $\mathcal{V}$ with the curvature $\Omega=
\frac{1}{2}R^{ab}\Sigma_{ab}\otimes\II_W
+\II_\Lambda\otimes\mathcal{F}\,,$ where $\mathcal{F}$ 
is the curvature of the Yang-Mills
connection.

In the following we will consider {\it homogeneous vector bundles} with
parallel bundle curvature, that is,
$
\nabla_\mu \mathcal{F}_{\alpha\beta}=0\,,
$
which means that the curvature satisfies the integrability constraints
\be 
[\mathcal{F}_{cd},\mathcal{F}_{ab}]  - R^f{}_{acd}\mathcal{F}_{fb} -
R^f{}_{bcd}\mathcal{F}_{af} = 0\,. 
\label{227} 
\ee


\subsection{Normal Coordinates}

Let $x'$ be a fixed point in $M$ and $\mathcal{U}$ be a sufficiently
small coordinate patch containing the point $x'$. Then every point $x$
in $\mathcal{U}$ can be connected with the point $x'$ by a unique
geodesic. We extend the local orthonormal frame $e_a{}^\mu(x')$ at the
point $x'$ to a local  orthonormal frame $e_a{}^\mu(x)$ at the point $x$
by parallel transport.  Of course, the frame $e_a{}^\mu$ depends on the
fixed point $x'$ as a parameter. Here and everywhere below the
coordinate indices of the tangent space at the point $x'$ are denoted by
primed Greek letters.  They are raised and lowered by the metric tensor
$g_{\mu'\nu'}(x')$ at the point $x'$.  The derivatives with respect to
$x'$ will be denoted by primed Greek indices as well. 

The parameters of the geodesic connecting the points $x$ and $x'$,
namely the unit tangent vector at the point $x'$  and the length of the
geodesic, (or, equivalently, the tangent vector at the point $x'$ with
the norm equal to the length of the geodesic), provide normal coordinate
system for $\mathcal{U}$. Let $d(x,x')$ be the geodesic distance between
the points $x$ and $x'$ and $\sigma(x,x')$ be a two-point function 
defined by $\sigma(x,x')=\frac{1}{2}[d(x,x')]^2$.
Then the derivatives $\sigma_{;\mu}(x,x')$ 
and $\sigma_{;\nu'}(x,x')$ are the tangent vectors to
the geodesic connecting the points $x$ and $x'$ at the points $x$
and $x'$ respectively pointing in opposite directions; one is obtained 
from another by parallel transport.
Now, let us define the quantities
$
y^a
=e^{a}{}_{\mu}\sigma^{;\mu}
=-e^{a}{}_{\mu'}\sigma^{;\mu'}\,.
$
These geometric parameters are nothing but the normal coordinates.


{\bf Remarks.}
Two remarks are in order here. First, strictly speaking, normal
coordinates can be only defined locally, in geodesic balls of radius
less than the injectivity radius of the manifold. However, for symmetric
spaces normal coordinates cover the whole manifold except for a set of
measure zero where they become singular \cite{camporesi90}. This set is
precisely the set of points conjugate to the fixed point $x'$
and of points that can be connected to the point
$x'$ by multiple geodesics. In any case, this set is a set of measure
zero and, as we will show below, it can be dealt with  by some
regularization technique. Thus, we will use the normal coordinates
defined above for the whole manifold. Second, for compact manifolds (or
for manifolds with compact submanifolds) the range of some normal
coordinates is also compact, so that if one allows them to range over
the whole real line $\RR$, then the corresponding compact submanifolds
will be covered infinitely many times.

\subsection{Curvature Group of a Symmetric Space}

We assumed that the manifold $M$ is locally symmetric. Since we also
assume that it is simply connected and complete, it is a globally
symmetric space (or simply symmetric space) \cite{wolf72}.  A symmetric
space is said to be compact, non-compact or Euclidean if all sectional
curvatures are positive, negative or zero. A generic symmetric space has
the structure $ M=M_0\times M_s\,, $ where $M_0=\RR^{n_0}$ and 
$M_s=M_+\times M_-$ is a semi-simple symmetric space; it is a product of
a compact symmetric space $M_+$ and a non-sompact symmetric space $M_-$.
Of course, the dimensions must satisfy the relation  $n_0+n_s=n$,
where $n_s=\dim M_s$.

The components of the curvature tensor
can be presented in the form \cite{avramidi96}
\be
R_{abcd} = \beta_{ik}E^i{}_{ab}E^k{}_{cd}\,,
\label{236}
\ee
where $E^i{}_{ab}$, $i=1,2,\dots,p$), is a collection
of $p$  anti-symmetric $n\times n$ matrices and
$\beta_{ik}$ is a symmetric
 nondegenerate  
$p\times p$ matrix with some $p\le n(n-1)/2$.

In the following the Latin indices from the middle of the alphabet will
be used to denote such matrices; they run over $1,\dots, p$ and
should not be
confused with the Latin indices from the beginning of the alphabet which
denote tensors in $M$. They will be raised and lowered with the matrix
$\beta_{ik}$ and its inverse $(\beta^{ik})$.

Next, we define the traceless $n\times n$
matrices $D_i=(D^a{}_{ib})$,
where
\be
D^a{}_{ib}=-\beta_{ik}E^k{}_{cb}\delta^{ca}\,.
\ee

The matrices $D_i$ are known to be the generators of the  holonomy
algebra, $\mathcal{H}$, i.e. the Lie algebra of the restricted holonomy 
group, $H$,
\be
[D_i, D_k] = F^j{}_{ik} D_j\,,
\label{310}
\ee
where $F^j{}_{ik}$ are the structure constants of the holonomy group.
The structure constants of the holonomy group define the
$p\times p$ matrices $F_i$, by $(F_i)^j{}_k=F^j{}_{ik}$, which generate
the adjoint representation  of the holonomy algebra,
\be
[F_i, F_k] = F^j{}_{ik} F_j\,.
\label{310xxx}
\ee

For symmetric spaces the introduced quantities satisfy additional
algebraic constraints. The most important consequence of the eq.
(\ref{312}) is the equation
\cite{avramidi96}
\be
E^i{}_{a c} D^c{}_{kb}
-E^i{}_{b c} D^c{}_{ka} 
= F^i{}_{kj}E^j{}_{ab}\,.
\label{313}
\ee
Now, by using the eqs. (\ref{310}) and (\ref{313}) one can prove 
that 
the matrix $\beta_{ik}$ satisfies the equation
\be
\beta_{ik}F^k{}_{jl}+\beta_{lk}F^k{}_{ji}=0\,.
\label{270}
\ee

Let $h^a{}_b$ be the projection to the subspace $T_x M_s$ of the tangent
space of dimension $n_s$, that is, the tensor $h_{ab}$ is nothing but
the metric tensor on the semi-simple subspace $T_xM_s$. Since the
curvature exists only in the semi-simple submanifold $M_s$, the
components of the curvature tensor $R_{abcd}$, as well as the tensors
$E^i{}_{ab}$,  are non-zero only in the semi-simple subspace $T_xM_s$.
Let $ q^a{}_{b}=\delta^a{}_b-h^a{}_b\, $ be the projection tensor to the
flat subspace $\RR^{n_0}$. Then 
\be
R_{abcd}q^a{}_e=R_{ab}q^a{}_e=E^i{}_{ab}q^a{}_e=D^a{}_{ib}q^b{}_e
=D^a{}_{ib} q_a{}^e=0\,.
\ee 

Now, we introduce a new type of indices, the capital Latin indices,
$A,B,C,\dots,$ which split according to $A=(a,i)$ and run from $1$ to
$N=p+n$. We define new quantities $C^A{}_{BC}$ by
\be
C^i{}_{ab}=E^i{}_{ab}, \qquad 
C^a{}_{ib}=-C^a{}_{bi}=D^a{}_{ib}, \qquad 
C^i{}_{kl}=F^i{}_{kl}\,,
\label{317}
\ee
all other components being zero. Then we can define $N\times N$ matrices
$C_A$ by  $(C_A)^B{}_C=C^B{}_{AC}$. 

Now, by
using the eqs. (\ref{310}),  (\ref{310xxx})
and (\ref{313}) one can prove the 
following theorem \cite{avramidi96}.

\begin{theorem}
The matrices $C_A$ satisfy the commutation
relations
\be
[C_A, C_B]=C^C{}_{AB}C_C\,,
\label{320}     
\ee
\end{theorem}

This means that the matricec $C_A$ generate the adjoint representation
of a Lie algebra $\mathcal{G}$  with the structure constants
$C^A{}_{BC}$. For the lack of a better name we call the algebra
$\mathcal{G}$ the {\it curvature algebra}. As it will be clear from the
next section it is a subalgebra of the total isometry algebra of the
symmetric space. It should be clear  that the holonomy algebra
$\mathcal{H}$ is the subalgebra of the curvature algebra  $\mathcal{G}$.  

Next, we define a symmetric nondegenerate
$N\times N$ matrix
\be
(\gamma_{AB}) = 
\left(
\begin{array}{cc}
\delta_{ab} & 0 \\
0 & \beta_{ik} \\
\end{array}
\right)
\,
\qquad
\mbox{and its inverse}
\qquad
(\gamma^{AB})=
\left(
\begin{array}{cc}
\delta^{ab} & 0 \\
0 & \beta^{ik} \\
\end{array}
\right)
\label{319}
\ee
These matrices 
will be
used to lower and to raise the capital Latin indices. 

Finally, by using the eqs. (\ref{313}) and (\ref{270}) one can show 
that
the matrix
$\gamma_{AB}$ satisfies the equation
\be
\gamma_{AB}C^B{}_{CD}+\gamma_{DB}C^B{}_{CA}=0\,.
\label{288xx}
\ee
Thus the curvature algebra $\mathcal{G}$ is compact; it is a direct sum
of two ideals, 
$
\mathcal{G}=\mathcal{G}_0\oplus \mathcal{G}_s,
$
an Abelian center $\mathcal{G}_0$ of dimension $n_0$ 
and a semi-simple algebra $\mathcal{G}_s$ of
dimension $p+n_s$.

\subsection{Killing Vectors Fields}

We will use extensively the isometries of the symmetric space $M$. 
The generators of isometries are the Killing vector fields. The set of
all Killing vector fields forms a representation of the  isometry
algebra, the Lie algebra of the isometry group of the manifold $M$. 
Let $K=(K^a{}_b)$ be the $n\times n$ matrix with the entries
\be
K^a{}_b=R^a{}_{cbd}y^cy^d.
\ee 
One can show \cite{avramidi96} that
\be
P_a=\left(\sqrt{K}\cot\sqrt{K}\right)^b{}_a\frac{\partial}{\partial y^b}\,,
\qquad
L_i=-D^b{}_{ia}y^a\frac{\partial}{\partial y^b}\,,
\label{229}
\ee
are Killing vector fields.
In the following we will only need these Killing vectors.
We introduce the following notation $(\xi_A)=(P_a,L_i)$. 

Next, by using the explicit
form of the Killing vector fields one can prove the following
theorem \cite{avramidi96}.
\begin{theorem}
The Killing vector fields $\xi_A$ satisfy the commutation relations
\be
[\xi_A, \xi_B]=C^C{}_{AB}\xi_C\,.
\label{320a}
\ee
\end{theorem}

Notice that they {\it do not} generate  the complete isometry algebra of
the symmetric space $M$ but rather they form a representation of the
curvature algebra $\mathcal{G}$ introduced in the previous section, which
is a subalgebra of the total isometry algebra.
It is clear that the Killing vector fields $L_i$ form a
representation of the holonomy algebra $\mathcal{H}$, which is the isotropy
algebra of the semi-simple submanifold $M_s$, and a  subalgebra of the
total isotropy algebra of the symmetric space $M$.

We list some properties of the Killing vector fields
that will be used below (for details, see \cite{avramidi07})
\be
\gamma^{AB}\xi_A{}^\mu\xi_B{}^\nu=g^{\mu\nu}\,,
\label{267}
\ee
\be
\gamma^{AB}\xi_A{}^\alpha\xi_B{}^\mu{}_{;\nu\lambda}
= R^\alpha{}_{\lambda\nu\mu}\,,
\label{235}
\ee
\be
\gamma^{AB}
\xi_A{}^\mu\xi_B{}^\nu{}_{;\beta}=0\,,
\label{2127a}
\ee
\be
\xi_A^c{}_{;a}\xi_B{}^b{}_{;c}
-\xi_B^c{}_{;a}\xi_A{}^b{}_{;c}
=C^C{}_{AB}\xi_C{}^b{}_{;a}-R^b{}_{acd}\xi_A{}^c\xi_B{}^d\,,
\label{2124x}
\ee
\be
\gamma^{AB}
\xi_A{}^\mu{}_{;\alpha}\xi_B{}^\nu{}_{;\beta}
=R^\mu{}_{\alpha}{}^{\nu}{}_{\beta}\,.
\label{242}
\ee

\subsection{Homogeneous Vector Bundles}

Equation (\ref{227}) imposes strong constraints on the curvature of the
homogeneous bundle $\mathcal{W}$. We define
\be
\mathcal{B}_{ab}=\mathcal{F}_{cd}q^c{}_a q^d{}_b\,,
\qquad
\mathcal{E}_{ab}=\mathcal{F}_{cd}h^c{}_bh^d{}_b\,,
\label{2139}
\ee
so that
$
\mathcal{B}_{ab}h^a{}_c=0\,,
$
and
$
\mathcal{E}_{ab}q^a{}_c=0\,.
$
Then, from eq. (\ref{227}) we obtain
\be
[\mathcal{B}_{ab}, \mathcal{B}_{cd}]
=[\mathcal{B}_{ab}, \mathcal{E}_{cd}]=0\,,
\ee
and
\be
[\mathcal{E}_{cd}, \mathcal{E}_{ab}] 
- R^f{}_{acd}\mathcal{E}_{fb}
- R^f{}_{bcd}\mathcal{E}_{af}
= 0\,.
\label{227a}
\ee
This means that  $\mathcal{B}_{ab}$ takes values in an Abelian ideal of the
gauge algebra $\mathcal{G}_{YM}$ and $\mathcal{E}_{ab}$ takes values in the
holonomy algebra. More precisely, eq. (\ref{227a}) is only possible if
the holonomy algebra $\mathcal{H}$ is an ideal of the gauge algebra 
$\mathcal{G}_{YM}$. Thus, the gauge group $G_{YM}$ must have a subgroup
$Z\times H$, where $Z$ is an Abelian group and $H$ is the holonomy
group.

The matrices $D^a{}_{ib}$ provide a
natural embedding of the holonomy algebra $\mathcal{H}$ in the orthogonal
algebra $\mathcal{SO}(n)$ in the following sense. Let $X_{ab}$ be the
generators of the orthogonal algebra $\mathcal{SO}(n)$ is some
representation.
Let $T_i$ be the matrices defined by
\be
T_i=-\frac{1}{2}D^a{}_{ib}X^b{}_a\,.
\ee
Then one can show that they
satisfy the commutation relations
\be
[T_i, T_k]=F^j{}_{ik}T_j\,.
\label{2127}
\ee
Thus $T_i$ are the generators of the gauge algebra $\mathcal{G}_{YM}$
realizing a representation $T$ of the holonomy algebra $\mathcal{H}$. 
Since $\mathcal{B}_{ab}$ takes values in the Abelian ideal of the algebra
of the gauge group we also have
\be
[\mathcal{B}_{ab},T_j]=0\,.
\ee

Then by using eq. (\ref{313}) one can show that
\footnote{We correct here a
sign misprint in eq. (3.24) in \cite{avramidi96}.}
\be
\mathcal{E}_{ab}=\frac{1}{2}R^{cd}{}_{ab}X_{cd}
=-E^i{}_{ab}T_i\,.
\ee
and
\bea
\mathcal{F}_{ab}&=&
-E^i{}_{ab}T_i+\mathcal{B}_{ab}
=\frac{1}{2}R^{cd}{}_{ab}X_{cd}+\mathcal{B}_{ab}
\,.
\label{2153mm}
\eea

Now, we consider the representation $\Sigma$ of the orthogonal algebra
defining the spin-tensor bundle $\mathcal{T}$ and define the matrices
\be
G_{ab}=\Sigma_{ab}\otimes\II_X+\II_\Sigma\otimes X_{ab}\,.
\ee
Obviously, these matrices are the generators of the orthogonal algebra
in the product representation $\Sigma\otimes X$.
Next, the matrices 
\be
Q_i=-\frac{1}{2}
D^a{}_{ib}\Sigma^b{}_a\,
\label{2104}
\ee
form a representation $Q$ of the holonomy algebra
$\mathcal{H}$,
and the matrices 
\bea
\mathcal{R}_{i}&=&
Q_i\otimes\II_{T}+\II_{\Sigma}\otimes T_i
=
-\frac{1}{2}D^a{}_{ib}G^b{}_a
\,
\label{2156mm}
\eea
are the generators of the holonomy algebra in the
product representation $\mathcal{R}=Q\otimes T$. 
Then the total curvature of  a
twisted spin-tensor bundle $\mathcal{V}$ is
\bea
\Omega_{ab}&=&-E^i{}_{ab}\mathcal{R}_i
+\mathcal{B}_{ab}
=\frac{1}{2}R^{cd}{}_{ab}G_{cd}+\mathcal{B}_{ab}
\,.
\label{2144}
\eea

\subsection{Twisted Lie Derivatives}

Let $\varphi$ be a section of a  twisted homogeneous spin-tensor bundle
$\mathcal{T}$. Let $\xi_A$ be the basis of Killing vector fields. Then the
covariant (or generalized, or twisted)  
Lie derivative of $\varphi$ along $\xi_A$ is defined by
\be
\mathcal{L}_A\varphi=
\mathcal{L}_{\xi_A}\varphi
=\left(\nabla_{\xi_A}+S_A\right)\varphi\,,
\label{2160x}
\ee
where $\nabla_{\xi_A}=\xi_A{}^\mu\nabla_\mu$, and
$
S_A=\frac{1}{2}\xi_A{}^a{}_{;b}G^b{}_a\,.
$
Note that
$
S_a q^a{}_b=0\,.
$

\begin{proposition}
There hold
\be
[\nabla_{\xi_A},\nabla_{\xi_B}]\varphi
=\left(
C^C{}_{AB}\nabla_{\xi_C}-\mathcal{R}_{AB}
+\mathcal{B}_{AB}
\right)\varphi\,,
\label{2160}
\ee
\be
\nabla_{\xi_A}S_B
=\mathcal{R}_{AB}\,,
\label{2152}
\ee
\be
[S_A,S_B]=C^C{}_{AB}S_C
-\mathcal{R}_{AB}\,,
\label{2165}
\ee
where
\bea
\mathcal{R}_{AB}&=&
\xi_A{}^a\xi_B{}^bE^i{}_{ab} \mathcal{R}_i
=
-\frac{1}{2}R^{cd}{}_{ab}\xi_A{}^a\xi_B{}^b G_{cd}
\,,
\eea
\be
\mathcal{B}_{AB}=\xi_A{}^a\xi_B{}^b\mathcal{B}_{ab}\,.
\ee
\end{proposition}

\noindent
\noindent{\bf Proof.}
By using the properties of the Killing vectors described in the previous
section and the eq. (\ref{2144}) we obtain first (\ref{2160}). Next,  we
obtain (\ref{2152}), and, further, by using the eq. (\ref{2124x}) we get
(\ref{2165}).

\rightline{$\Box$}

We define the  operator
\be
\mathcal{L}^2=\gamma^{AB}\mathcal{L}_A\mathcal{L}_B\,.
\ee

\begin{theorem}
The operators $\mathcal{L}_A$ and $\mathcal{L}^2$
satisfy the commutation relations
\be
[\mathcal{L}_A, \mathcal{L}_B] 
= C^C{}_{AB} \mathcal{L}_C+\mathcal{B}_{AB},
\label{2173xx}
\ee
\be
[\mathcal{L}_A,\mathcal{L}^2]
=2\gamma^{BC}\mathcal{B}_{AB}\mathcal{L}_C\,.
\label{2178xx}
\ee
\end{theorem}

\noindent
\noindent\noindent\noindent\noindent{\bf Proof.} 
This follows from 
\be
[\mathcal{L}_A,\mathcal{L}_B]=[\nabla_{\xi_A},\nabla_{\xi_B}]
+[\nabla_{\xi_A},S_B]
-[\nabla_{\xi_B},S_A]
+[S_A,S_B]
\ee
and eqs. (\ref{2160}), (\ref{2152}), and (\ref{2165}). The eq.
(\ref{2178xx}) follows directly from (\ref{2173xx}).

The operators $\mathcal{L}_A$ form an algebra that
is a direct sum of a nilpotent ideal and a 
semisimple algebra. For the lack of a better name we call this algebra
{\it gauged curvature algebra} and denote it by $\mathcal{G}_{\rm gauge}$.

Now, by using the eqs. (\ref{2127a}),
(\ref{2152}) and (\ref{235}) one can prove that
\be
\gamma^{AB}\xi_A{}^\mu S_B=0\,,
\qquad
\gamma^{AB}\nabla_{\xi_A} S_{B}=0\,,
\qquad
\gamma^{AB}S_A S_B=\mathcal{R}^2\,.
\label{244xxx}
\ee

We define the Casimir operator
\be
{\cal R}^2=\beta^{ij}{\cal R}_i{\cal R}_j
=\frac{1}{4}R^{abcd}G_{ab}G_{cd}\,.
\ee

\begin{theorem}
The Laplacian $\Delta$ acting on sections of a twisted spin-tensor
bundle $\mathcal{V}$ over a symmetric space has the form
\be
\Delta=\mathcal{L}^2-\mathcal{R}^2\,.
\label{2110}
\ee
Therefore,
\be
[\mathcal{L}_{A},\Delta] = 2 \gamma^{BC}\mathcal{B}_{AB}\mathcal{L}_C\,.
\label{43}
\ee
\end{theorem}

\noindent
\noindent\noindent\noindent\noindent{\bf Proof.} We have
\be
\gamma^{AB}\mathcal{L}_A\mathcal{L}_B=\gamma^{AB}\nabla_{\xi_A}\nabla_{\xi_B}
+\gamma^{AB}S_A\nabla_{\xi_B}
+\gamma^{AB}\nabla_{\xi_A}S_B
+\gamma^{AB}S_AS_B\,.
\ee
Now, by using eqs. (\ref{267}) and (\ref{2127a}) we get
\be
\gamma^{AB}\nabla_{\xi_A}\nabla_{\xi_B}=\Delta\,.
\ee
Next, by using the eqs. (\ref{2152}) and (\ref{244xxx}) we obtain
(\ref{2110}). The eq. (\ref{43}) follows from the commutation relations 
(\ref{2173xx}).

\subsection{Isometries and Pullbacks}

Let $\omega^i$ be the canonical coordinates on the holonomy group and 
$(k^A)=(p^a,\omega^i)$ be the canonical coordinates on the  gauged
curvature group.  
Let
$\xi=\left<k,\xi\right>=k^A\xi_A
=p^a P_a+\omega^iL_i$ be a Killing vector field
and let $\psi_t: M\to M$ be the
one-parameter diffeomorphism (the isometry) generated by the vector
field  $\xi$. Let $\hat x=\psi_t(x)$, so that
$
\frac{d\hat x}{dt}=\xi{}(\hat x)\,
$
and
$
\hat x\big|_{t=0}=x\,.
$
The solution of this equation 
$
\hat x=\hat x(t,p,\omega,x,x')\,.
$
depends on the parameters 
$t,p,\omega,x$ and $x'$.
We will be interested mainly in the case when the points $x$ and $x'$
are close to each other. In fact, at the end of our calculations
we will take the limit $x=x'$.

Now, we choose the normal coordinates $y^a$ of the point   defined above
and the normal coordinates $\hat y^a$ of the point $\hat x$  with the
origin at $x'$, so that  the normal coordinates $y'$ of the point $x'$
are  equal to zero. Recall that the normal coordinates are
equal to the  components of the tangent vector at the point $x'$  to the
geodesic connecting the points $x'$ and the current point, that is,
$y^a=-e^a{}_{\mu'}(x')\sigma^{;\mu'}(x,x')$ and $\hat
y^a=-e^a{}_{\mu'}(x')\sigma^{;\mu'}(\hat x,x')$. Then by taking into
account the explicit form of the Killing vectors given by eq.
(\ref{229}) we have
\be
{d \hat y^a\over dt} = \left(\sqrt {K(\hat y)}\cot \sqrt
{K(\hat y)}\right)^a_{\ b}p^b -\omega^iD^a{}_{ib}\hat y^b\,,
\label{447}
\ee
with the initial condition
$
\hat y^a\big|_{t=0}=y^a\,.
$
The solution of this equation defines a function 
$\hat y=\hat y(t,p,\omega,y)$. 

We define the matrix
\be
D(\omega)=\omega^i D_i\,.
\ee

\begin{proposition}
The Taylor expansion of the function $\hat y=\hat y(t,p,\omega,y)$
in $p$ and $y$ reads
\be
\hat y^a = \left(\exp[-tD(\omega)]\right)^a{}_b y^b
+\left({1-\exp[-t D(\omega)]\over D(\omega)}\right)^a{}_bp^b
+O(y^2,p^2,py)\,.
\label{2197mm}
\ee
There holds
\be
\det\left(
\frac{\partial \hat y^a}{\partial p^b}
\right)\Bigg|_{p=y=0,t=1}
=\det{}_{TM}\left({\sinh[\,D(\omega)/2]\over D(\omega)/2}\right)\,.
\label{429a}
\ee
\end{proposition}
\noindent{\bf Proof.}
First, for $p=0$ from eq. (\ref{447}) 
we obtain
$
\hat y(t,0,\omega,0)=0\,.
$
Next, by differentiating the eq. (\ref{447}) with respect to 
$y^b$ and setting $p=y=0$ we obtain
\be
\frac{\partial \hat y^a}{\partial y^b}\Bigg|_{p=y=0}
=\left(\exp[-tD(\omega)]\right)^a{}_b\,.
\label{2204mm}
\ee

Further, by differentiating the eq. (\ref{447}) with respect to $p^b$
and setting $p=0$, we obtain a differential equation
for the matrix
$
\frac{\partial \hat y^a}{\partial p^b}\Big|_{p=y=0}
$
whose 
solution is
\be
\left(\frac{\partial \hat y^a}{\partial p^b}\right)\Big|_{p=y=0}
={1-\exp[-t D(\omega)]\over D(\omega)}\,.
\label{452}
\ee
By using the obtained results
we get the desired
formula (\ref{2197mm}).
Finally, 
by taking into account that the matrix $D(\omega)$ is traceless,
by using eq. (\ref{452}) we obtain (\ref{429a}).

\rightline{$\Box$}

The function $\hat y=\hat y(t,p,\omega,y)$
implicitly defines the function
$
p=p(t,\omega,\hat y,y)\,.
$
The function 
$\bar p=\bar p(\omega,y)$ is now defined
by the equation
$
\hat y(1,\bar p,\omega,y)=0\,,
$
or
$
\bar p(\omega,y)=p(1,\omega,0,y)\,.
$

\begin{proposition}
The Taylor expansion of the function $\bar p(\omega,y)$
in $y$ has the form
\be
\bar p^a=-\left(D(\omega)
\frac{\exp[-D(\omega)]}{1-\exp[-D(\omega)]}
\right)^a{}_b y^b
+O(y^2)\,.
\label{2213mm}
\ee
Therefore,
\be
\det\left(-
\frac{\partial \bar p^a}{\partial y^b}
\right)\Bigg|_{y=0}
=\det{}_{TM}\left({\sinh[\,D(\omega)/2]
\over D(\omega)/2}\right)^{-1}\,.
\label{429aa}
\ee
\end{proposition}
\noindent{\bf Proof.}
Next, by taking into account the eq. 
$\hat y\big|_{p=y=0}=0$
 we have
$
\bar p\Big|_{y=0}=0\,.
$
Further, by differentiating the equation
$\hat y(1,\bar p,\omega,y)=0$
with respect to
$y^c$ and setting $y=0$ we get
\be
\frac{\partial \hat y^a}{\partial y^b}\Big|_{p=y=0, t=1}
+\frac{\partial \hat y^a}{\partial p^c}\Big|_{p=y=0, t=1}
\frac{\partial \bar p^c}{\partial y^b}\Big|_{y=0}=0\,,
\ee
and, therefore,
\be
\frac{\partial \bar p^a}{\partial y^b}\Big|_{y=0}
=-\left(D(\omega)
\frac{\exp[-D(\omega)]}{1-\exp[-D(\omega)]}
\right)^a{}_b\,.
\ee
This leads to both (\ref{2213mm}) and (\ref{429aa}).

\rightline{$\Box$}

Now, we define
\be
\Lambda^{\hat\mu}{}_{\nu}=
\frac{\partial \hat x^\mu}{\partial x^{\nu}}\,,
\qquad
\mbox{so that}
\qquad
(\Lambda^{-1})^{\mu}{}_{\hat\alpha}
=g^{\mu\nu}(x)\Lambda^{\hat\beta}{}_{\nu}
g_{\hat\beta\hat\alpha}(\hat x)\,.
\ee
Let $e^a{}_\mu$ and $e_a{}^\mu$ be a local orthonormal frame 
that is obtained by parallel transport along geodesics from a
point $x'$ and $\psi_t^*$ be the pullback of the isometry
$\psi_t$ defined above.
Then the frames 
of $1$-forms $e^a$ and $\psi_t^* e^a$ are related by an 
orthogonal transformation
$
(\psi_t^* e^a)(x)=O^a{}_b e^b(x)\,,
$
where the matrix $O^a{}_b$ is defined by
\be
O^a{}_b
=e^a{}_{\hat\alpha}(\hat x)
\Lambda^{\hat\alpha}{}_{\mu} e_b{}^{\mu}(x)
\,.
\label{2149xxz}
\ee
Since the matrix $O$ is orthogonal, it can be
parametrized by
$
O=\exp\theta\,,
$
where $\theta_{ab}$ is an antisymmetric matrix.

\begin{proposition}
For $p=y=0$ the matrix $O$ has the form
\be
O\Big|_{p=y=0}=\exp\left[-t D(\omega)\right]\,.
\label{2149}
\ee
\end{proposition}
\noindent\noindent\noindent\noindent{\bf Proof.}
We use normal coordinates $\hat y^a$ and $y^a$. Then the matrix
$O$ takes the form
\be
O^a{}_b=e^a{}_{\hat\alpha}
\frac{\partial \hat x^\alpha}{\partial \hat y^c}
\frac{\partial \hat y^c}{\partial y^d}
\frac{\partial y^d}{\partial x^\mu}e_b{}^\mu\,.
\ee
Now, by using the explicit form of the Jacobian
in normal coordinates
and the fact that $\hat y\big|_{p=y=0}=0$
we obtain
\be
\frac{\partial y^a}{\partial x^\mu}e_b{}^\mu
\Bigg|_{p=y=0}
=e^a{}_{\hat\alpha}
\frac{\partial \hat x^\alpha}{\partial \hat y^b}
\Bigg|_{p=y=0}
=\delta^a{}_b\,.
\ee
Therefore,
\be
O^a{}_b\Big|_{p=y=0}
=\frac{\partial \hat y^a}{\partial y^b}\Bigg|_{p=y=0}\,,
\ee
and, finally (\ref{2204mm}) gives the desired result
(\ref{2149}).

\rightline{$\Box$}

Let $\varphi$ be a section of the twisted spin-tensor bundle $\mathcal{V}$.
Let $V_{x}$ be the fiber at the point $x$ and
$V_{\hat x}$ be the fiber at the point $\hat x=\psi_t(x)$.
The pullback of the diffeomorphism $\psi_t$ defines the map,
that we call just the pullback,
$
\psi_t^*: C^\infty(\mathcal{V})\to C^\infty(\mathcal{V})
$
on smooth sections of the twisted spin-tensor bundle $\mathcal{V}$.

\begin{proposition}
Let $\varphi$ be a section of a twisted spin-tensor bundle
$\mathcal{V}$. Then
\be
(\psi_t^*\varphi)(x)
=\exp\left(-\frac{1}{2}\theta_{ab}G^{ab}\right)
\varphi(\hat x)\,.
\label{2236mm}
\ee
In particular, for $p=y=0$ (or $x=x'$)
\be
(\psi_t^*\varphi)(x)\Big|_{p=y=0}
=\exp\left[t\mathcal{R}(\omega)\right]
\varphi(x')\,,
\label{2237mm}
\ee
where
$
\mathcal{R}(\omega)=\omega^i\mathcal{R}_i\,.
$
\end{proposition}
\noindent{\bf Proof.}
First, from the eq. (\ref{2149}) we see that
$
\theta^a{}_b\Big|_{p=y=0}=-t\omega^i D^a{}_{ib}\,.
$
Then, from the definition (\ref{2156mm})
of the matrices $\mathcal{R}_i$ we get (\ref{2237mm}).

\rightline{$\Box$}

\section{Heat Semigroup}
\setcounter{equation}0

\subsection{Geometry of the Curvature Group}

Let $G_{\rm gauge}$ be the gauged curvature group and $H$ be its
holonomy subgroup.  Both these groups have compact algebras. However,
while the holonomy group is always compact, the curvature group is, in
general, a product of a nilpotent group, $G_0$, and a
semi-simple group, $G_s$,
$
G_{\rm gauge}=G_0\times G_s\,.
$
The semi-simple group $G_s$ is a product $G_s=G_+\times G_-$ of a
compact $G_+$ and a non-compact $G_-$ subgroups.

Let $\xi_A$ be the basis Killing vectors, $k^A$ be the canonical
coordinates on the curvature group $G$ and $\xi(k)=k^A\xi_A$.
The canonical coordinates are exactly the normal coordinates
on the group defined above.
Let $C_A$ be the generators of the curvature group in
adjoint representation and $C(k)=k^AC_A$. 
In the following $\partial_M$ means the partial derivative 
$\partial/\partial k^M$ with
respect to the canonical coordinates.
We define the matrix $Y^A{}_M$ by the equation
\be
\exp[-\xi(k)]\partial_M\exp[\xi(k)]=Y^A{}_M\xi_A\,,
\label{31xx}
\ee
which is well defined since the right hand side lies in the Lie
algebra of the curvature group.
The matrix $Y=(Y^A{}_M)$ can be computed
explicitly, namely,
\be
Y=\frac{1-\exp[-C(k)]}{C(k)}\,.
\ee
Let $X=(X_A{}^M)=Y^{-1}$ be the inverse matrix of $Y$.
Then we define the $1$-forms $Y^A$ and the vector fields $X_A$
on the group $G$ by
\be
Y^A=Y^A{}_M dk^M\,,\qquad
X_A=X_A{}^M\partial_M\,.
\ee

\begin{proposition}
There holds
\be
X_A \exp[\xi(k)]=\exp[\xi(k)]\xi_A\,.
\label{35xx}
\ee
\end{proposition}
\noindent{\bf Proof.} This follows immediately from the eq. (\ref{31xx}).

\rightline{$\Box$}

Next, 
by differentiating the eq. (\ref{31xx}) 
with respect to $k^L$ and alternating the indices $L$ and $M$ we obtain
\be
\partial_L Y^A{}_M-\partial_M Y^A{}_L=-C^A{}_{BC}Y^B{}_{L}Y^C{}_M\,,
\label{36xx}
\ee
which, of course, can also be written as
\be
dY^A=-\frac{1}{2}C^A{}_{BC}Y^B\wedge Y^C\,.
\ee

\begin{proposition}
The vector fields $X_A$ satisfy the commutation relations
\be
[X_A,X_B]=C^C{}_{AB}X_C\,.
\label{37xx}
\ee
\end{proposition}
\noindent{\bf Proof.} This follows from the eq. (\ref{36xx}).

The vector fields $X_A$ are nothing but the right-invariant
vector fields. They form a representation  of
the curvature algebra.

We will also need the following fundamental property of Lie groups.

\begin{proposition}
Let $G$ be a Lie group with the structure constants $C^A{}_{BC}$,
$C_A=(C^B{}_{AC})$ and $C(k)=C_Ak^A$. Let $\gamma=(\gamma_{AB})$ be a
symmetric non-degenerate matrix satisfying the  equation
\be
(C_A)^T=-\gamma C_A \gamma^{-1}\,.
\ee
Let $X=(X_A{}^M)$ be a matrix defined by
\be
X=\frac{C(k)}{1-\exp[-C(k)]}\,.
\ee
Then
\be
(\det X)^{-1/2}\gamma^{AB}X_A{}^M\partial_M X_B{}^N\partial_N
(\det X)^{1/2} = -\frac{1}{24}\gamma^{AB}C^C{}_{AD}C^D{}_{BC}\,.
\label{425}
\ee
\end{proposition}
\noindent{\bf Proof.} It is easy to check that this equation holds at $k=0$. Now,
it can be proved by showing that it is a group invariant. For a 
detailed proof for semisimple groups see
\cite{helgason84,camporesi90,dowker71}.

\rightline{$\Box$}

It is worth stressing that this equation holds not only on semisimple
Lie groups but on any group with a compact Lie algebra, that is, when
the structure constants $C^A{}_{BC}$ and the matrix $\gamma_{AB}$, used
to define the metric $G_{MN}$ and the operator $X^2$, satisfy the eq.
(\ref{288xx}). Such algebras can have an Abelian center.

Now, by using the right-invariant 
vector fields we define a metric on 
the curvature group $G$
\be
G_{MN} = \gamma_{AB} Y^{A}{}_{M} Y^{B}{}_{N}\,,
\qquad
G^{MN} = \gamma^{AB} X_A{}^{M} X_B{}^{N}\,.
\label{418}
\ee
This metric is bi-invariant.
This means that the vector fields $X_A$ are the Killing vector fields
of the metric $G_{MN}$.
One can easily show that 
this metric defines the following natural affine
connection $\nabla^G$ on the group
\be
\nabla^G_{X_C} X_A=-\frac{1}{2}C^A{}_{BC}X_B{}\,,
\qquad
\nabla^G_{X_C} Y^A{}=\frac{1}{2}C^B{}_{AC}Y^B{}\,,
\label{311xx}
\ee
with the scalar curvature
\be
R_G=-\frac{1}{4}\gamma^{AB}C^C{}_{AD}C^D{}_{BC}\,.
\label{312xx}
\ee

Since the matrix $C(k)$ is traceless we have $\det\exp[C(k)/2]=1$,
and, therefore, the volume element on the group is
\be
|G|^{1/2}=\left(\det G_{MN}\right)^{1/2}
=|\gamma|^{1/2}\det{}_{\mathcal{G}}
\left({\sinh[C(k)/2]\over C(k)/2}\right)\,,
\label{419}
\ee
where $|\gamma|=\det \gamma_{AB}$.

It is not difficult to see that
\be
k^M Y^A{}_M=k^M X_M{}^A=k^A\,.
\label{314xx}
\ee
By differentiating this equation with respect to $k^B$ and 
contracting the indices $A$ and $B$ we obtain
\be
k^M\partial_A X_M{}^A=N-X_A{}^A\,.
\label{315xx}
\ee

Now, by contracting the eq. (\ref{311xx}) with $G^{BC}$ we obtain
the zero-divergence condition for the right-invariant vector fields
\be
|G|^{-1/2}\partial_M\left(|G|^{1/2}X_A{}^M\right)=0\,.
\label{313xx}
\ee
Next, we define the Casimir operator
\be
X^2=C_2(G,X)=\gamma^{AB}X_AX_B\,.
\ee
By using the eq. (\ref{313xx}) one can easily show that $X^2$ is
an invariant differential operator that is nothing
but the scalar Laplacian on the group
\be
X^2=|G|^{-1/2}\partial_M |G|^{1/2}G^{MN}\partial_N
=G^{MN}\nabla^G_M\nabla^G_N\,.
\ee
Then, by using the eqs. (\ref{288xx}) and (\ref{320}) 
one can show that the operator $X^2$ commutes with the operators
$X_A$,
\be
[X_A,X^2]=0\,.
\ee

Since we will actually be working with the gauged curvature group,
we introduce now the operators (covariant right-invariant
vector fields) $J_A$ by
\be
J_A=X_A-\frac{1}{2}\mathcal{B}_{AB}k^B\,,
\ee
and the operator 
\be
J^2=\gamma^{AB}J_AJ_B\,.
\ee

\begin{proposition}
The operators $J_A$ and $J^2$ satisfy the commutation relations
\be
[J_A,J_B]=C^C{}_{AB}J_C+\mathcal{B}_{AB}\,,
\label{326xx}
\ee
and
\be
[J_A,J^2]=2\mathcal{B}_{AB}J^B\,.
\label{326xxx}
\ee
\end{proposition}
\noindent{\bf Proof.}
By using the eq. 
(\ref{2139})
we obtain
\be
X_B{}^A \mathcal{B}_{AM}
=
\gamma_{BN}\gamma^{AC}X_C{}^N \mathcal{B}_{AM}
=\mathcal{B}_{BM}\,,
\label{327xx}
\ee
and, hence,
\be
\gamma^{AB}X_B{}^M \mathcal{B}_{AM}
=0\,,
\ee
and, further, by using (\ref{37xx})
we obtain (\ref{326xx}).
By using the eqs. (\ref{327xx}) we get (\ref{326xxx}).

Thus, the operators $J_A$ form a representation of the gauged curvature
algebra.
Now, let $\mathcal{L}_A$ be the operators of Lie derivatives
satisfying the commutation relations (\ref{2173xx}) and
$\mathcal{L}(k)=k^A \mathcal{L}_A$. 

\begin{proposition}
There holds
\be
J_A \exp[\mathcal{L}(k)]=\exp[\mathcal{L}(k)]\mathcal{L}_A\,.
\ee
and, therefore,
\be
J^2 \exp[\mathcal{L}(k)]=\exp[\mathcal{L}(k)]\mathcal{L}^2\,.
\label{330xx}
\ee
\end{proposition}
\noindent{\bf Proof.}
We have
\be
\exp[-\mathcal{L}(k)]\partial_M\exp[\mathcal{L}(k)]
=\exp[-Ad_{\mathcal{L}(k)}]\partial_M\,.
\ee
By using the commutation relations (\ref{2173xx}) and eq. (\ref{2139})
we obtain
\be
\exp[-\mathcal{L}(k)]\partial_M\exp[\mathcal{L}(k)]
=Y^A{}_M\mathcal{L}_A+\frac{1}{2}\mathcal{B}_{MN}k^N\,.
\ee
The statement of the proposition follows from the definition of the
operators $J_A$, $J^2$ and $\mathcal{L}^2$.

\subsection{Heat Kernel on the Curvature Group}

Let $\mathcal{B}$ be the matrix with the components 
$\mathcal{B}=(\gamma^{AB}\mathcal{B}_{BC})$ so that
\be
\mathcal{B}=(\mathcal{B}_{AB})=\left(
\begin{array}{cc}
\mathcal{B}_{ab} & 0 \\
0 & 0\\
\end{array}
\right)
\label{2172}
\ee
Let $k^A$ be the canonical coordinates on the curvature group $G$
and $A(t;k)$ be a function defined by
\be
A(t;k)=\det{}_\mathcal{G}
\left(\frac{\sinh\left[C(k)/2+t\mathcal{B}\right]}
{C(k)/2+t\mathcal{B}}\right)^{-1/2}\,.
\ee
By using the eqs. (\ref{327xx}) one can rewrite this in the form
\be
A(t;k)=\det{}_\mathcal{G}
\left(\frac{\sinh\left[C(k)/2\right]}
{C(k)/2}\right)^{-1/2}
\det{}_\mathcal{G}
\left(\frac{\sinh\left[t\mathcal{B}\right]}
{t\mathcal{B}}\right)^{-1/2}
\,.
\label{334xx}
\ee
Notice also that
\be
\det{}_\mathcal{G}
\left(\frac{\sinh\left[t\mathcal{B}\right]}
{t\mathcal{B}}\right)^{-1/2}
=\det{}_{TM}
\left(\frac{\sinh\left[t\mathcal{B}\right]}
{t\mathcal{B}}\right)^{-1/2}\,,
\ee
where $\mathcal{B}$ is now regarded as just the matrix 
$\mathcal{B}=(\mathcal{B}^a{}_b)$.

Let $\Theta(t;k)$ be another function on the group $G$ defined by
\be
\Theta(t;k)=\frac{1}{2}\left<k, \gamma\hat\Theta k\right>\,,
\label{337mm}
\ee
where $\hat\Theta$ is the matrix
\be
\hat\Theta=t\mathcal{B}\coth(t\mathcal{B})\,
\label{337xx}
\ee
and $\left<u,\gamma v\right>=\gamma_{AB}u^Av^B$ is the inner product on
the algebra $\mathcal{G}$.

\begin{theorem}
Let $\Phi(t;k)$ be a function on the group $G$ defined by
\bea
\Phi(t;k)&=&(4\pi t)^{-N/2}A(t;k)
\exp\left(-\frac{\Theta(t;k)}{2t}+\frac{1}{6}R_G t\right)\,,
\label{46a}
\eea
Then $\Phi(t;k)$ satisfies the equation
\be
\partial_t \Phi = J^2\Phi\,,
\label{423}
\ee
and the initial condition
\be
\Phi(0;k)=|\gamma|^{-1/2}\delta(k)\,.
\label{428}
\ee
\end{theorem}

\noindent
\noindent\noindent\noindent\noindent{\bf Proof.}
We compute first
\be
\partial_t\Phi=
\left[
\frac{1}{6}R_G-\frac{1}{2t}\tr_\mathcal{G}\hat\Theta
+\frac{1}{4t^2}\left<k,\gamma\hat\Theta^2k\right>
-\frac{1}{4}\left<k,\gamma\mathcal{B}^2k\right>
\right]\Phi\,.
\label{339xx}
\ee
Next, we have
\be
J^2=X^2-\gamma^{AB}\mathcal{B}_{AC}k^C X_B
+\frac{1}{4}\gamma^{AB}\mathcal{B}_{AC}\mathcal{B}_{BD}k^Ck^D\,.
\ee
By using the eqs. (\ref{327xx}) and  (\ref{2172})
and the anti-symmetry of the
matrix $\mathcal{B}_{AB}$ we show that
\be
\mathcal{B}_{AC}k^C X_B\Phi=0\,.
\ee
Thus,
\bea
J^2\Phi&=&
\Biggl[
A^{-1}(X^2A)
-\frac{1}{2t}(X^2\Theta)
+\frac{1}{4t^2}\gamma^{AB}(X_A\Theta)(X_B\Theta)
\nonumber\\
&&
-\frac{1}{t}A^{-1}\gamma^{AB}(X_B A) (X_A \Theta)
-\frac{1}{4}\left<k,\gamma\mathcal{B}^2k\right>
\Biggr]\Phi\,.
\label{344xx}
\eea

Further, by using \ref{327xx}) we get
\be
\gamma^{AB}(X_A\Theta)(X_B\Theta)
=\left<k,\gamma\hat\Theta^2 k\right>\,,
\qquad
X^2\Theta=\tr_\mathcal{G}X+\tr_\mathcal{G}\hat\Theta-N\,.
\label{345xx}
\ee
Now, by using the eq. (\ref{313xx})
and eqs. (\ref{2172}) and (\ref{315xx}) we show that
\be
A^{-1}\gamma^{AB}(X_A\Theta)X_B A=\frac{1}{2}\left(
N-\tr_\mathcal{G}X\right)\,,
\ee
and by using eq. (\ref{425}) we obtain
\be
A^{-1}X^2 A =\frac{1}{6}R_G\,.
\label{349xx}
\ee

Finally, substituting the eqs. (\ref{345xx})-(\ref{349xx}) into eq.
(\ref{344xx}) and comparing it with eq. (\ref{339xx}) we prove the
eq. (\ref{423}). The initial condition (\ref{428}) follows easily from
the well known property of the Gaussian.
This completes the proof of the theorem.

\subsection{Regularization and Analytical Continuation}

In the following we will complexify the gauged  curvature group in the
following sense. We extend the  canonical coordinates
$(k^A)=(p^a,\omega^i)$ to the whole complex Euclidean space  $\CC^N$.
Then all group-theoretic functions introduced above become analytic
functions of $k^A$ possibly with some poles on the real section $\RR^N$
for compact groups.  In fact, we replace the actual real slice $\RR^N$
of $\CC^N$ with an $N$-dimensional subspace $\RR^N_{\rm reg}$ in $\CC^N$
obtained by rotating the real section $\RR^N$ counterclockwise in
$\CC^N$ by $\pi/4$. That is, we replace each coordinate $k^A$ by
$e^{i\pi/4}k^A$. In the complex domain the group becomes non-compact. We
call this procedure the {\it decompactification}. If the group is
compact, or has a compact subgroup, then this plane will cover the
original group infinitely many times. 

Since the metric $(\gamma_{AB})=\diag(\delta_{ab},\beta_{ij})$ is not
necessarily positive definite,  (actually, only the metric of the
holonomy group $\beta_{ij}$ is non-definite) we analytically continue
the function $\Phi(t;k)$ in the  complex plane of $t$ with a cut along
the negative imaginary axis so that
$-\pi/2<\arg\,t<3\pi/2$.  Thus, the function $\Phi(t;k)$ defines an
analytic function of $t$ and $k^A$. For the purpose of the following
exposition we shall consider $t$ to be {\it real negative}, $t<0$. This
is needed in order to make all integrals convergent and well defined and
to be able to do the analytical continuation.

As we will show below, the singularities occur only in the holonomy
group. This means that there is no need to complexify the coordinates
$p^a$. Thus, in the following we assume the coordinates $p^a$ to be real
and the coordinates $\omega^i$ to be complex, more precisely, to take
values in the $p$-dimensional subspace $\RR^p_{\rm reg}$ of $\CC^p$ 
obtained by rotating $\RR^p$ counterclockwise by $\pi/4$ in $\CC^p$ That
is, we have $\RR^N_{\rm reg}=\RR^n\times \RR^p_{\rm reg}$.

This procedure (that we call a regularization) with the nonstandard
contour of integration is necessary  for the convergence of the
integrals below since we are treating both the compact and the
non-compact symmetric spaces simultaneously. Recall, that, in general,
the nondegenerate diagonal matrix $\beta_{ij}$ is not positive definite.
The space $\RR^p_{\rm reg}$ is chosen in such a way to make the Gaussian
exponent purely imaginary. Then the indefiniteness of the matrix $\beta$
does not cause any problems. Moreover, the integrand does not have any
singularities on these contours. The convergence of the integral is
guaranteed by the exponential growth of the sine for imaginary argument.
These integrals can be computed then in the following way. The
coordinates $\omega^j$ corresponding to the compact directions  are
rotated further by another $\pi/4$ to imaginary axis and the coordinates
$\omega^j$ corresponding to the non-compact directions are rotated back
to the real axis. Then, for $t<0$ all the integrals below are well
defined and convergent and define an analytic function of $t$ in a complex
plane with a cut along the negative imaginary axis.

\subsection{Heat Semigroup}

\begin{theorem}
The heat semigroup $\exp(t\mathcal{L}^2)$
can be represented 
in form of the integral
\be
\exp(t\mathcal{L}^2) = 
\int\limits_{\RR^N_{\rm reg}} dk\; |G|^{1/2}(k)\Phi(t;k)
\exp[\mathcal{L}(k)]\,.
\label{49a}
\ee
\end{theorem}

\noindent\noindent\noindent\noindent{\bf Proof.} Let
\be
\Psi(t) = 
\int\limits_{\RR^N_{\rm reg}} dk\; |G|^{1/2}\Phi(t;k)
\exp[\mathcal{L}(k)]\,.
\label{490}
\ee
By using the previous theorem we obtain
\be
\partial_t\Psi(t)
=\int\limits_{\RR^N_{\rm reg}} dk\; |G|^{1/2}\exp[\mathcal{L}(k)]
J^2\Phi(t;k)\,.
\label{421}
\ee
Now, by integrating by parts we get
\be
\partial_t\Psi(t)
=\int\limits_{\RR^N_{\rm reg}} dk\; |G|^{1/2}\Phi(t;k)
J^2\exp[\mathcal{L}(k)]\,,
\label{421a}
\ee
and, by using eq. (\ref{330xx}) we obtain
\be
\partial_t\Psi(t)=\Psi(t)\mathcal{L}^2\,.
\ee
Finally from the initial condition (\ref{428}) for the function
$\Phi(t;k)$ we get
$
\Psi(0)=1\,,
$
and, therefore, $\Psi(t)=\exp(t\mathcal{L}^2)$.

\begin{theorem}
Let $\Delta$ be the Laplacian acting on sections of a homogeneous
twisted spin-tensor vector bundle over a symmetric space. Then the
heat semigroup $\exp(t\Delta)$
can be represented in form of an integral
\bea
\exp(t\Delta) 
&=& 
(4\pi t)^{-N/2}
\det{}_{TM}\left(\frac{\sinh(t\mathcal{B})}{t\mathcal{B}}\right)^{-1/2} 
\exp\left(-t\mathcal{R}^2+ {1\over 6} R_G t\right)
\nonumber\\
&&
\times\int\limits_{\RR^N_{\rm reg}} dk\; |\gamma|^{1/2}
\det{}_\mathcal{G}\left({\sinh[C(k)/2]\over C(k)/2}\right)^{1/2}
\nonumber\\
&&
\times
\exp\left\{ -{1\over 4t}
\left<k,\gamma t\mathcal{B}\coth(t\mathcal{B})k\right>
\right\}
\exp[\mathcal{L}(k)]\,.
\label{49}
\eea
\end{theorem}

\noindent
\noindent\noindent\noindent{\bf Proof.}
By using the eq. (\ref{2110}) we obtain
\be
\exp(t\Delta)=\exp\left(-t\mathcal{R}^2\right)
\exp\left(t\mathcal{L}^2\right)\,.
\ee
The statement of the theorem follows now from the eqs. (\ref{49a}),
(\ref{46a}),  (\ref{334xx})-(\ref{337xx}) and (\ref{419}).

\section{Heat Kernel}
\setcounter{equation}0

\subsection{Heat Kernel Diagonal and Heat Trace}

The heat kernel diagonal on a homogeneous bundle over a symmetric space
is parallel. In a parallel local frame it is just a constant 
matrix. The fiber trace of the heat kernel diagonal is just a constant.
That is why, it can be computed at any point in $M$. We fix a point $x'$
in $M$ and choose the normal coordinates $y^a$ with the origin at $x'$
such that the Killing vectors 
are given by the explicit formulas above
(\ref{229}). We compute the heat kernel diagonal at the
point $x'$.

The heat kernel diagonal can be obtained by acting by the heat
semigroup $\exp(t\Delta)$ on the delta-function,
\cite{avramidi94,avramidi96}
\bea
U^{\rm diag}(t)&=&\exp(t\Delta)\delta(x,x')\Big|_{x=x'}
\nonumber\\
&=&
\exp\left(-t\mathcal{R}^2\right)
\int\limits_{\RR^N_{\rm reg}} dk\; |G|^{1/2}\Phi(t;k)
\exp[\mathcal{L}(k)]\delta(x,x')\Big|_{x=x'}\,.
\label{41}
\eea
To be able to use this integral representation we need to compute the
action of the isometries $\exp[\mathcal{L}(k)]$ on the delta-function.

It is not very difficult to check that
the Lie derivatives are nothing but 
the generators of the pullback, that is,
\be
\mathcal{L}_{\xi}\varphi
=k^A\mathcal{L}_A\varphi
=\frac{d}{dt}(\psi_t^*\varphi)
\Big|_{t=0}\,.
\ee
By using this fundamental fact we can prove the following
proposition.

\begin{proposition}
Let $\varphi$ be a section of the twisted spin-tensor bundle 
$\mathcal{V}$, $\mathcal{L}_A$ be the twisted Lie derivatives,
$k^A=(p^a,\omega^i)$ be the canonical coordinates on the group
and
$\mathcal{L}(k)=k^A\mathcal{L}_A$.
Let $\xi=k^A\xi_A$ be the Killing vector and $\psi_t$
be the corresponding one-parameter diffeomorphism.
Then
\bea
\exp\left[\mathcal{L}(k)\right]\varphi(x)
&=&\exp\left(-\frac{1}{2}\theta_{ab}G^{ab}\right)
\varphi(\hat x)\Bigg|_{t=1}\,,
\label{2135}
\eea
where $\hat x=\psi_t(x)$ and the matrix $\theta$ is defined above 
by $O=\exp\theta$ where $O$ is defined by (\ref{2149xxz}).
In particular, for $p=0$ and $x=x'$
\bea
\exp[\mathcal{L}(k)]\varphi(x)\Big|_{p=0,x=x'}
&=&
\exp\left[\mathcal{R}(\omega)\right]
\varphi(x)\,.
\label{2150}
\eea
\end{proposition}
\noindent\noindent\noindent\noindent{\bf Proof.}
This statement follows from eqs. (\ref{2236mm}) and (\ref{2237mm}).

\begin{proposition}
Let $\omega^i$ be the canonical coordinates on the holonomy
group $H$ and $(k^A)=(p^a,\omega^i)$ be the natural splitting of the
canonical coordinates on the curvature group $G$.
Then
\be
\exp[\mathcal{L}(k)]\delta(x,x')\Big|_{x=x'}
=\det{}_{TM}\left(\frac{\sinh[D(\omega)/2]}{D(\omega)/2}\right)^{-1}
\exp[\mathcal{R}(\omega)]
\delta(p)\,.
\label{424a}
\ee
\end{proposition}

\noindent
\noindent\noindent\noindent\noindent{\bf Proof.}
Let $\hat x(t,p,\omega,x,x')=\psi_t(x)$.
By making use of the eq. (\ref{2135})
we obtain
\bea
\exp[\mathcal{L}(k)]\delta(x,x')\Big|_{x=x'}
&=&\exp\left(-\frac{1}{2}\theta_{ab}G^{ab}\right)
\delta(\hat x(1,p,\omega,x,x'),x')\Big|_{x=x',t=1}\,.
\eea
Now we change the variables from $x^\mu$ 
to the normal coordinates $y^a$
to get
\be
\delta(\hat x(1,p,\omega,x,x'),x')\Big|_{x=x'}
=|g|^{-1/2}
\det\left(\frac{\partial y^a}{\partial x^\mu}\right)
\delta(\hat y(1,p,\omega,y))\Big|_{y=0}\,.
\ee
This delta-function picks the values of $p$ that make
$\hat y=0$, which is exactly the functions
$\bar p=\bar p(\omega,y)$
defined above
by $\hat y(1,\bar p,\omega,y)=0$. By switching further
to the  variables $p$ we obtain
\be
\delta(\hat x(1,p,\omega,x,x'),x')\Big|_{x=x'}
=|g|^{-1/2}
\det\left(\frac{\partial y^a}{\partial x^\mu}\right)
\det\left(\frac{\partial \hat y^b}{\partial p^c}\right)^{-1}
\delta(p-\bar p(\omega,y))\Big|_{y=0,t=1}\,.
\ee
Now, by recalling that 
$\bar p|_{y=0}=0$ and by using the Jacobian 
$\partial y^a/\partial x^\nu$
of the transformation 
to the normal coordinates in symmetric spaces
and (\ref{429a})
we evaluate the Jacobians for $p=y=0$ and $t=1$ to 
get the eq. (\ref{424a}).

{\bf Remarks.}
Some remarks are in order here. 
We implicitly assumed that there are no closed geodesics
and that the equation of closed orbits of isometries
$
\hat y^a(1,\bar p,\omega,0)=0\,
$
has a unique solution $\bar p=\bar p(\omega,0)=0$. On compact symmetric
spaces this is not true: there are infinitely many closed geodesics and
infinitely many closed orbits of isometries.   However, these global
solutions, which reflect the global topological structure of the
manifold, will not affect our local analysis.  In particular, they do
not affect the asymptotics of the heat kernel.  That is why, we have
neglected them here.  This is reflected in the fact that the Jacobian in
(\ref{424a}) can become singular when the coordinates of the holonomy
group $\omega^i$ vary from $-\infty$ to $\infty$.  Note that the exact
results for compact symmetric spaces can be obtained by an analytic
continuation from the dual noncompact case when such closed geodesics
are absent \cite{camporesi90}. That is why we proposed above to
complexify our holonomy group. If the coordinates $\omega^i$ are complex
taking values in the subspace $\RR^p_{\rm reg}$ defined above, then the
equation $\hat y^a(1,\bar p,\omega,0)=0$
should have a unique solution and the Jacobian is
an analytic function. It is worth stressing once again  that the
canonical coordinates cover the whole group except for a set of measure
zero.  Also a compact subgroup is covered  infinitely many times. 

Now by using the above lemmas and the theorem we can compute the heat
kernel diagonal. 
We define the
matrix
\be
F(\omega)=\omega^i F_i\,,
\ee
and
\be
R_H=-\frac{1}{4}\beta^{ik}F^j{}_{im}F^m{}_{kj}\,,
\label{295vv}
\ee
which is nothing but the scalar curvature of the isotropy group $H$.

\begin{theorem}
The heat kernel diagonal of the  Laplacian on twisted spin-vector
bundles over a symmetric space has the form
\bea
U^{\rm diag}(t)
&=&(4\pi t)^{-n/2}
\det{}_{TM}\left(\frac{\sinh(t\mathcal{B})}{t\mathcal{B}}\right)^{-1/2}
\exp\left\{\left({1\over 8} R + {1\over 6} R_H 
-\mathcal{R}^2\right)t\right\}
\nonumber
\label{437a}
\\
&&
\times
\int\limits_{\RR^n_{\rm reg}} 
\frac{d\omega}{(4\pi t)^{p/2}}\;|\beta|^{1/2}
\exp\left\{-{1\over 4 t}\left<\omega,\beta\omega\right>\right\}
\cosh\left[\,\mathcal{R}(\omega)\right]
\nonumber\\
&& 
\times\det{}_\mathcal{H}
\left({\sinh\left[\,F(\omega)/2\right]\over 
\,F(\omega)/2}\right)^{1/2}
\det{}_{TM}\left({\sinh\left[\,D(\omega)/2\right]\over 
\,D(\omega)/2}\right)^{-1/2}\,,
\eea
where $|\beta|=\det \beta_{ij}$ and
$\left<\omega,\beta\omega\right>=\beta_{ij}\omega^i\omega^j$.
\end{theorem}

\noindent
\noindent\noindent\noindent\noindent{\bf Proof.}
First, we have $dk=dp\;d\omega$ and
$|\gamma|=|\beta|\,.$
By using the equations (\ref{41}) 
and (\ref{424a}) and integrating over $p$ we obtain 
the heat kernel diagonal
\bea
U^{\rm diag}(t)
&=&\int\limits_{\RR^p_{\rm reg}}d\omega\;|G|^{1/2}(0,\omega)
\Phi(t;0,\omega)
\det{}_{TM}\left(\frac{\sinh[D(\omega)/2]}{D(\omega)/2}\right)^{-1}
\nonumber\\
&&\times
\exp[\mathcal{R}(\omega)-t\mathcal{R}^2]\,.
\eea
Further, by using the definition of the matrices $C_A$
we compute the determinants
\be
\det{}_\mathcal{G}\left({\sinh[C(\omega)/2]\over C(\omega)/2}\right)
=\det{}_{TM}\left({\sinh[D(\omega)/2]\over D(\omega)/2}\right)
\det{}_\mathcal{H}\left({\sinh[F(\omega)/2]\over F(\omega)/2}\right)\,.
\label{437aa}
\ee
Now, by using (\ref{2172}) we compute (\ref{337mm}) 
$
\Theta(t;0,\omega)=\frac{1}{2}\left<\omega,\beta\omega\right>\,,
$
and, finally, 
by using eq. (\ref{46a}), (\ref{334xx}), (\ref{312xx})
and (\ref{295vv}) 
we get the result (\ref{437a}).

\rightline{$\Box$}

By using this theorem we can also compute the heat trace
for compact symmetric spaces 
\bea
\Tr_{L^2}\exp(t\Delta)
=\int\limits_M d\vol\; \tr_V U^{\rm diag}(t)
&=& \vol(M)
\tr_V U^{\rm diag}(t)
\label{437cc}
\nonumber
\eea
where $\tr_V$ is the fiber trace.

\subsection{Heat Kernel Asymptotics}

It is well known that there is the following asymptotic
expansion as $t\to 0$ of the heat kernel diagonal
\cite{gilkey95}
\be
U^{\rm diag}(t)\sim (4\pi t)^{-n/2}
\sum_{k=0}^\infty t^k a_k\,.
\ee
The coefficients $a_k$ are called the local heat kernel
coefficients.
On compact symmetric spaces there is a similar
asymptotic expansion of the heat trace
\be
\Tr_{L^2}\exp(t\Delta)\sim (4\pi t)^{-n/2}
\sum_{k=0}^\infty t^k A_k\,
\ee
with the global heat invariants $A_k$ defined by
\be
A_k=\int_M d\vol\; \tr_V a_k=\vol(M)\tr_V a_k\,.
\ee

We introduce a Gaussian average over the holonomy algebra by
\be
\left<f(\omega)\right> = \int\limits_{\RR^p_{\rm reg}}
\frac{d\omega}{(4\pi)^{p/2}}\; |\beta|^{1/2}
\exp\left(-{1\over 4}\left<\omega,\beta\omega\right>\right)f(\omega)
\ee
Then we can write
\bea
&&U^{\rm diag}(t) 
=(4\pi t)^{-n/2}  
\det{}_{TM}\left(\frac{\sinh(t\mathcal{B})}{t\mathcal{B}}\right)^{-1/2}
\exp\left\{\left({1\over 8} R + {1\over 6} R_H 
-\mathcal{R}^2\right) t \right\}
\label{422}
\\[12pt]
&&
\times
\Bigg< \cosh\left[\sqrt{t}\,\mathcal{R}(\omega)\right]
\det{}_\mathcal{H}\left({\sinh\left[\sqrt{t}\,F(\omega)/2\right]\over 
\sqrt{t}\,F(\omega)/2}\right)^{1/2}
\det{}_{TM}\left({\sinh\left[\sqrt{t}\,D(\omega)/2\right]\over 
\sqrt{t}\,D(\omega)/2}\right)^{-1/2}\Bigg>   
\nonumber
\eea
This equation can be used now to generate all heat kernel
coefficients $a_k$ for any locally symmetric space  simply by expanding
it in a power  series in $t$. By using the standard Gaussian averages
\bea
\left<\omega^i_1\cdots \omega^{i_{2k+1}}\right> &=& 0\,,
\label{440a}
\\
\left<\omega^{i_1}\cdots \omega^{i_{2k}}\right> 
&=& {(2k)!\over k!}\beta^{(i_1 i_2}\cdots
\beta^{i_{2k-1}i_{2k})}\,,
\label{441a}
\eea
one can obtain now all heat kernel coefficients in terms of 
traces of various contractions 
of the matrices $D^a{}_{ib}$ and $F^j{}_{ik}$ with  the
matrix  $\beta^{ik}$.
All these quantities are curvature invariants and
can be  expressed directly in terms of the Riemann tensor.

\section{Conclusion}

We have continued the study of the heat kernel on homogeneous spaces
initiated in 
\cite{avramidi93,avramidi94,avramidi95,avramidi96}. In those
papers we have developed a systematic technique for calculation of the
heat kernel in two cases: a) a Laplacian on a vector bundle  with a
parallel curvature over a flat space \cite{avramidi93,avramidi95}, and
b) a scalar Laplacian on manifolds with parallel curvature
\cite{avramidi94,avramidi96}. What was missing in that study was the
case of  a non-scalar Laplacian on vector bundles with parallel
curvature over curved manifolds with parallel curvature. 

In the present paper we considered the Laplacian on a homogeneous bundle
and generalized the technique developed in \cite{avramidi96} to compute
the corresponding heat semigroup and the heat kernel. It is worth
pointing out that our formal result applies to general symmetric spaces
by making use of the regularization and the analytical continuation
procedure described above. Of course, the heat kernel coefficients are
just polynomials in the curvature and do not depend on this kind of
analytical continuation (for more details, see \cite{avramidi96}).  

As we mentioned above, due to existence of  multiple closed geodesics
the obtained form of the heat kernel for compact symmetric spaces 
requires an additional regularization, which consists simply in an
analytical continuation of the result from the complexified noncompact
case. In any case, it gives a generating function for all heat
invariants and reproduces correctly the whole asymptotic expansion of
the heat kernel diagonal. However, since there are no closed geodesics
on non-compact symmetric spaces, it seems that the analytical
continuation of the obtained result for the heat kernel diagonal should
give the {\it exact result} for the non-compact case, and, even more
generally, for the general case too.

\section*{Acknowledgements}

I would like to thank the organizers of the special 2007 Midwest
Geometry Conference in the honor of Thomas P. Branson, in particular, 
Palle Jorgensen, for the invitation and for the financial support. I
would like to dedicate this contribution to the memory of  Thomas P.
Branson (1953-2006) whose sudden death shocked all his friends,
collaborators and relatives. Tom was a good friend and a  great
mathematician who had a major impact on many areas of modern
mathematics.


\end{document}